\documentclass[12pt,twoside,leqno]{article}
\begin{document}

\author{Constantin Udri\c{s}te and Teodor Oprea}
\title{H-convex Riemannian submanifolds }
\date{}
\maketitle

\begin{abstract}
Having in mind the well known model of Euclidean convex hypersurfaces [4],
[5], and the ideas in [1] many authors defined and investigate convex
hypersurfaces of a Riemannian manifold. As it was proved by the first author
in [7], there follows the interdependence between convexity and Gauss
curvature of the hypersurface. In this paper we define $H$-$convexity$ of a
Riemannian submanifold of arbitrary codimension, replacing the normal versor
of a hypersurface with the mean curvature vector. A
characterization, used by B.Y. Chen [2], [3] as the definition of \textsl{%
strictly }$H$-$convexity$, it is obtained.
\end{abstract}

\section{Convex hypersurfaces in Riemannian ma\-nifolds}

Let $(N,g)$ be a complete finite-dimensional Riemannian manifold and $M$ be
an oriented hypersurface whose induced Riemannian is also denoted by $g.$
Let $x$ be a point in $M\subset N$ and $V$ a neighborhood of $x$ in $N$ such
that $\exp _{x}:T_{x}N\rightarrow V$ is a diffeomorphism.

We denote by $\omega $ the $1$-form associated to the unit normal vector
field $\xi $ on $M$. The real-valued function defined on $V$ by 
\[
F(y)=\omega _{x}(\exp _{x}^{-1}(y)) 
\]
has the property that the set 
\[
TGH_{x}=\{y\in V\left| F(y)=0\right. \} 
\]
is a \textit{totally geodesic hypersurface at }$x$, tangent to\textit{\ }$M%
\mathit{\ }$at\textit{\ }$x$.

This hypersurface is the common boundary of the sets

$$TGH_{x}^{-}=\{y\in V\left|F(y)\leq 0\}, \right.
TGH_{x}^{+}=\{y\in V\left| F(y)\geq 0\}.\right.
$$

\textbf{Definition. }The hypersurface $M$ is called \textsl{convex }at%
\textsl{\ }$x\in M$ if there exists an open set $U\subset V\subset N$
containing $x$ such that $M\cap U$ is contained either in $TGH_{x}^{-},$ or
in $TGH_{x}^{+}$.

A hypersurface $M$ convex at $x$ is said to be \textsl{strictly convex} at $%
x $ if

\[
M\cap U\cap TGH_{x}=\{x\}. 
\]

Let us recall some results in [6]-[7]. We start with a necessary
condition for a hypersurface of a Riemannian manifold to be convex at a
given point.

Let $M$ be an oriented hypersurface in the manifold $N$, let $\xi$ be the normal versor and $h$ be the second fundamental form of $M$.

\textbf{Theorem 1.1} (Udri\c{s}te)\textbf{.} {\it If $M$ is convex at $x\in M,$ then the bilinear form 

$$
\Omega _{x}:T_{x}M\times T_{x}M\rightarrow R, 
\Omega _{x}(X,Y)=g(h(X,Y),\xi ), 
$$
is semidefinite.}

\textsl{\ \medskip }

The converse of Theorem 1.1 is not true. To show this, we consider the
surface $M:x^{3}=(x^{1})^{2}+(x^{2})^{3}$ in $R^{3}$. One observes that $%
0\in M$, $\xi (0)=(0,0,1)$ and $TGH_{0}:$ $x^{3}=0$ is the plane tangent to $%
M$ at the origin. On the other hand, if 
$$
c:I\rightarrow M, 
c(t)=(x^{1}(t),x^{2}(t),x^{3}(t)), 
$$
is a $C^{2}$ curve such that $0\in I$ and $c(0)=0$, then $(x^{3})^{\prime
\prime }(0)=2((x^{1})^{\prime }(0))^{2}$ and hence the function 
$$
f:I\rightarrow R, 
f(t)=\left\langle c(t)-0,\xi (0)\right\rangle 
$$
satisfies the relations $f(t)=x^{3}(t)$ and $f^{\prime \prime
}(0)=(x^{3})^{\prime \prime }(0)=2((x^{1})^{\prime }(0))^{2}.$

Since $\Omega _{0}(\stackrel{.}{c}(0),\stackrel{.}{c}(0))=f^{\prime \prime
}(0)$, and $c$ is an $C^{2}$ arbitrary curve, one gets that $\Omega _{0}$ is
positive semidefinite.

\medskip

However $M$ is not convex at the origin because the tangent plane\newline
$TGH_{0}:$ $x^{3}=0$ cuts the surface along the semicubic parabola 
\[
x^{3}=0,(x^{1})^{2}+(x^{2})^{3}=0 
\]
and consequently in any neighborhood of the origin there exist points of the
surface placed both below the tangent plane and above the tangent
plane.\medskip

\medskip

If the bilinear form $\Omega $ is definite at the point $x\in M$, then the
hypersurface $M$ is strictly convex at $x$.

\medskip

The next results establish a connection between the Riemannian manifolds
admitting a function whose Hessian is positive definite and their convex
hypersurfaces.

\medskip

\textbf{Theorem 1.2} (Udri\c{s}te). \textsl{Suppose that the Riemannian
manifold }$(N,g)$\textsl{\ supports a function }$f:N\rightarrow R$\textsl{\
with positive definite Hessian. On each compact oriented hypersurface }$M$%
\textsl{\ in }$N$\textsl{\ there exists a point }$x\in M$\textsl{\ such that
the bilinear form }$\Omega (x)$\textsl{\ is definite.}

\medskip

\textbf{Theorem 1.3} (Udri\c{s}te). \textsl{\ Suppose that the Riemannian
manifold }$(N,g)$\textsl{\ supports a function }$f:N\rightarrow R$\textsl{\
with positive definite Hessian. Then}

\textsl{1) There is no compact minimal hypersurface in }$N$\textsl{.}

\textsl{2) If the hypersurface }$M$\textsl{\ is connected and compact and
its Gauss curvature is nowhere zero, then }$M$\textsl{\ is strictly convex.}

\medskip

\textbf{Theorem 1.4} (Udri\c{s}te). \textsl{Let }$(N,g)$\textsl{\ be a
connected and complete Riemannian manifold and }$f:N\rightarrow R$\textsl{\
a function with positive definite Hessian. If }$x_{0}$\textsl{\ is a
critical point of }$f$\textsl{\ and a}$_{0}=f(x_{0})$\textsl{, then for any
real number }$a\in $\textsl{Im}$f\backslash $\textsl{\ }$\{a_{0}\}$\textsl{\
the hypersurface }$M_{a}=f^{-1}\{a\}$\textsl{\ is strictly convex}.

\section{H-convex Riemannian submanifolds}

Having in mind the model of convex hypersurfaces in Riemannian manifolds,
we define $H-convexity$ of a Riemannian submanifold of arbitrary
codimension, replacing the normal versor $\xi$ of a hypersurface with the mean
curvature vector $H$.

\medskip

Let $(N,g)$ be a complete finite-dimensional Riemannian manifold and $M$ be
a submanifold in $N$ of dimension $n$ whose induced Riemannian is also
denoted by $g.$ We denote by $H$ the mean curvature vector field of $M$ and by $\omega $ the associated $1$-form. 

Let $x$ be a point in $M\subset N$, with $H_{x}\neq 0$ and $V$ a
neighborhood of $x$ in $N$ such that $\exp _{x}:T_{x}N\rightarrow V$ is a
diffeomorphism.
The real-valued function defined on $V$ by 
\[
F(y)=\omega _{x}(\exp _{x}^{-1}(y)) 
\]
has the property that the set 
\[
TGH_{x}=\{y\in V\left| F(y)=0\right. \} 
\]
is a \textit{totally geodesic hypersurface at }$x$, tangent to\textit{\ }$M%
\mathit{\ }$at\textit{\ }$x$.

This hypersurface is the common boundary of the sets

$$
TGH_{x}^{-}=\{y\in V\left| F(y)\leq 0\},\right. \\ 
TGH_{x}^{+}=\{y\in V\left| F(y)\geq 0\}.\right.
$$

\textbf{Definition. }The submanifold $M$ is called \textsl{H-convex }at%
\textsl{\ }$x\in M$ if there exists an open set $U\subset V\subset N$
containing $x$ such that $M\cap U$ is contained either in $TGH_{x}^{-},$ or
in $TGH_{x}^{+}$.

A submanifold $M,$ $H$-convex at $x$ is said to be \textsl{strictly }$H$-%
\textsl{convex} at $x$ if

\[
M\cap U\cap TGH_{x}=\{x\}. 
\]

The next result is a necessary condition for a submanifold of a Riemannian
manifold to be $H$-convex at a given point.

\medskip

\textbf{Theorem 2.1.}\textsl{\ If }$M$\textsl{\ is a submanifold in }$N,$ 
\textsl{H-convex at} $x\in M,$ \textsl{then the bilinear form} 
$$
\Omega _{x}:T_{x}M\times T_{x}M\rightarrow R, 
\Omega _{x}(X,Y)=g(h(X,Y),H), 
$$
\textsl{where }$h$\textsl{\ is the second fundamental form of }$M,$\textsl{\
is positive semidefinite.}

\textsl{\medskip }

\textbf{Proof. }We suppose that there is a open set $U\subset V\subset N$
which contains $x$ such that $M\cap U$ $\subset $ $TGH_{x}^{+}$.
For an arbitrary vector $X\in T_{x}M$ let $c:I\rightarrow M\cap U$ be a $%
C^{2}$ curve, where $I$ is a real interval such that $0\in I$ and $c(0)=x$, $%
\stackrel{.}{c}(0)=X.$

As $c(I)\subset M\cap U$ $\subset $ $TGH_{x}^{+}$ the function $f=F\circ
c:I\rightarrow R$ satisfies\newline
(1) $f(t)\geq 0$, $\forall $ $t\in I$.\ 

It follows that $0$ is a global minimum point for $f$ , therefore\newline
(2) $0=f^{^{\prime }}(0)=\omega _{x}(d\exp _{x}^{-1}(c(0)))(\stackrel{.}{c}%
(0))=\omega _{x}(X)$,\newline
(3) $0\leq f^{^{\prime \prime }}(0)=\omega _{x}(d^{2}\exp _{x}^{-1}(c(0)))(%
\stackrel{.}{c}(0),\stackrel{.}{c}(0))+$\newline
$+\omega _{x}(d\exp _{x}^{-1}(c(0)))(\stackrel{..}{c}(0))=\omega _{x}(%
\stackrel{..}{c}(0))=\Omega _{x}(X,X).$

Since $X\in T_{x}M$ is an arbitrary vector, we obtain that $\Omega _{x}$ is
positive semidefinite.

\textbf{Remark.} We consider $\{e_{1},e_{2},...,e_{n}\}$ an orthonormal
frame in $T_{x}M$. Since Trace($\Omega
_{x})=g(\sum\limits_{i=1}h(e_{i},e_{i}),H_{x})=ng(H_{x},H_{x})>0$, the quadratic form $%
\Omega _{x}$ cannot be negative semidefinite, therefore $M\cap U$ cannot be
contained in $TGH_{x}^{-}.$ So, if the submanifold $M$ is $H$-convex at the
point $x$, then there exists an open set $U\subset V\subset N$ containing $x$
such that $M\cap U$ is contained in $TGH_{x}^{+}$.

In the sequel, we intend to prove that if the bilinear form $\Omega _{x}$ is
positive definite, then the submanifold $M$ is strictly $H$-convex at the
point $x$. For this purpose we introduce a function similar to the height function
used in the study of the hypersurfaces of an Euclidean space.
We fix $x\in M\subset N$ and a neighborhood $V$ of $x$ for which $\exp
_{x}:T_{x}N\rightarrow V$ is a diffeomorphism. The function 
$$
F_{\omega _{x}}:V\rightarrow R,
F_{\omega _{x}}(y)=\omega _{x}(\exp _{x}^{-1}(y) )
$$
has the property that it is affine on geodesics radiating from $x$.

We consider an arbitrary vector $X\in T_{x}M$ and a curve $c:I\rightarrow V$
such that $0\in I$, $c(0)=x$, $\stackrel{.}{c}(0)=X$. 
The function $f=F_{\omega _{x}}\circ c:I\rightarrow R$ satisfy
$$f^{\prime }(0)=\omega _{x}(d\exp _{x}^{-1}(c(0))(\stackrel{.}{c}%
(0))=\omega _{x}(\stackrel{.}{c}(0))=\omega _{x}(X)=g(H,X)=0,$$
and hence $x\in M$ is a critical point of $F_{\omega _{x}}$.

\textbf{Theorem 2.2.} \textsl{Let }$M$\textsl{\ be a submanifold in }$N$%
\textsl{. If the bilinear form }$\Omega _{x}$\textsl{\ is positive definite,
then }$M$\textsl{\ is strictly }$H$\textsl{-convex at the point }$x$\textsl{.%
}

\medskip

\textbf{Proof.} The point $x\in M$ is a critical point of $F_{\omega _{x}}$
and $F_{\omega _{x}}(x)=0.$ On the other hand one observe that 
$$
{Hess}^{N}F_{\omega _{x}}={Hess}^{M}F_{\omega _{x}}-dF_{\omega _{x}}(\Omega H). 
$$
As $F_{\omega _{x}}$ is affine on each geodesic radiating from $x$, it
follows Hess$^{N}F_{\omega _{x}}=0.$ It remains that 
\[
{Hess}^{M}F_{\omega _{x}}(x)=\Omega _{x} 
\]
and hence Hess$^{M}F_{\omega _{x}}$ is positive definite at the point $x$.
In this way $x$ is a strict local minimum point for $F_{\omega _{x}}$ in $%
M\cap V$, i. e., the submanifold $M$ is strictly $H$-convex at $x$.

\medskip

\textbf{Remark.} 1) The bilinear form $\Omega _{x}$ is positive
(semi)definite if and only if the Weingarten operator $A_{H}$ is positive
(semi)definite.

2) If $M$ is a hypersurface in $N$, $x$ is a point in $M$ with $H_{x}\neq 0$%
, then $M$ is $H$-convex at $x$ if and only if $M$ is convex at $x$.

A class of strictly $H$-convex submanifolds into a Riemannian manifold is
made by the curves which have the mean curvature nonzero.

\medskip

\textbf{Theorem 2.3.} \textsl{Let }$(N,g)$\textsl{\ be a Riemannian manifold
and }$c:I\rightarrow N$\textsl{\ a regular curve which have the mean
curvature nonzero, where }$I$\textsl{\ is a real interval. Then }$c$\textsl{%
\ is a strictly }$H$\textsl{-convex submanifold of }$N$\textsl{.}

\medskip

\textbf{Proof.} We fix $t\in I$. As $T_{c(t)}c=$Sp$\{\stackrel{.}{c}(t)\}$,
we obtain 
\[
H_{c(t)}=\frac{h(\stackrel{.}{c}(t),\stackrel{.}{c}(t))}{\left\| \stackrel{.%
}{c}(t)\right\| ^{2}},
\]
hence $\Omega (\stackrel{.}{c}(t),\stackrel{.}{c}(t))=g(h(\stackrel{.}{c}(t),%
\stackrel{.}{c}(t)),H_{c(t)})=\left\| \stackrel{.}{c}(t)\right\| ^{2}\left\|
H_{c(t)}\right\| ^{2}>0$, therefore $\Omega $ is positive definite. It
follows that $c$\textsl{\ }is a strictly $H$-convex submanifold of $N$.

\section{H-convex Riemannian submanifolds in real space forms}

Let us consider $(M,g)$ a Riemannian manifold of dimension $n$. We fix $x\in
M$ and $k\in \overline{2,n}.$ Let $L$ be a vector subspace of dimension $k$ in $T_{x}M$. If $X\in L$ is a
unit vector, and $\{e_{1}^{\prime },e_{2}^{\prime },...,e_{k}^{\prime }\}$
an orthonormal frame in $L$, with $e_{1}^{\prime }=X$, we shall denote

\[
{Ric}_{L}(X)=\sum\limits_{j=2}^{k}k(e_{1}^{\prime }\wedge
e_{j}^{\prime }), 
\]
where $k(e_{1}^{\prime }\wedge e_{j}^{\prime })$ is the sectional curvature
given by Sp$\{e_{1}^{\prime },e_{j}^{\prime }\}.$ 
We define \textsl{the Ricci curvature of k-order} at the point $x\in M$,

\[
\theta _{k}(x)=\frac{1}{k-1}\min\limits _{L,\dim L=k, \atop \\ X\in L,
\left\| X\right\| =1}  {Ric}_{L}(X), 
\]

B.Y. Chen showed in [2], [3] that the eigenvalues of the Weingarten operator
of a submanifold in a real space form and the Ricci curvature of $k$-order
satisfies the inequality in the next 

\textbf{Theorem 3.1} (Chen) \textsl{Let }$(\widetilde{M}(c),\widetilde{g})$%
\textsl{\ be a real space form of dimension }$m$\textsl{\ and }$M\subset 
\widetilde{M}(c)$\textsl{\ a submanifold of dimension }$n$\textsl{, and }$%
k\in \overline{2,n}$\textsl{. Then}

\textsl{\ }

i)\textsl{\ }$A_{H}\geq \frac{n-1}{n}(\theta _{k}(x)-c)I_{n}.$

\textsl{\ }

ii)\textsl{\ If }$\theta _{k}(x)\neq c,$\textsl{\ then the previous
inequality is strict. }

\medskip

\textbf{Corollary.} \textsl{If }$M$\textsl{\ is a submanifold of dimension }$%
n$\textsl{\ in the real space form }$\widetilde{M}(c)$\textsl{\ of dimension 
}$m$\textsl{, }$x\in M$\textsl{\ and there is a natural number }$k\in 
\overline{2,n}$\textsl{\ such that }$\theta _{k}(x)>c$\textsl{, then }$M$%
\textsl{\ is strictly} $H$-\textsl{convex at the point }$x$.

The converse of previous corollary is also true in the case of hypersurfaces
in a real space form.

\textbf{Theorem 3.2.} \textsl{If }$M$\textsl{\ is a hypersurface of
dimension }$n$\textsl{\ of a real space form }$\widetilde{M}(c),$\textsl{\ }$%
x\in M$\textsl{\ and }$M$\textsl{\ is strictly }$H$-\textsl{convex at the
point }$x$\textsl{, then} 
\[
\theta _{k}(x)>c,\forall \ k\in \overline{2,n}.
\]

\medskip

\textbf{Proof.} Let $x$ be a point in $M,$ $H$ the mean curvature of $M$ and 
$\pi $ a $2$-plane in $T_{x}M$. We consider $\{X,Y\}$ an orthonormal frame in 
$\pi $ and $\xi =\frac{H_{x}}{\left\| H_{x}\right\| }.$
The second fundamental form of the submanifold $M$ satisfies the relation%
$$
h(U,V)=\frac{\Omega _{x}(U,V)}{\left\| H_{x}\right\| }\,\,\xi , \forall \,\,U,V\in T_{x}M.\leqno (1)
$$

From the Gauss equation, one gets

$$
\widetilde{R}(X,Y,X,Y)=R(X,Y,X,Y)-\widetilde{g}(h(X,X),h(Y,Y))+
+\widetilde{g}(h(X,Y),h(X,Y)).\leqno (2)
$$

Using the relation (1) and the fact that $\widetilde{M}(c)$ has the
sectional curvature $c$, we obtain
$$
R(X,Y,X,Y)=c+\frac{1}{\left\| H_{x}\right\| ^{2}}(\Omega
_{x}(X,X)\Omega _{x}(Y,Y)-\Omega _{x}(X,Y)^{2}).\leqno (3)
$$

On the other hand $\Omega _{x}$ is positive definite because $M$ is strictly 
$H$-convex at the point $x$. From the Cauchy inequality, using the fact that $X$
and $Y$ are linear independents, we obtain 
$$
\Omega _{x}(X,X)\Omega _{x}(Y,Y)-\Omega (X,Y)^{2}>0.\leqno (5)
$$

From (3) and (4) we find
(5) $R(X,Y,X,Y)>c,$ which means that the sectional curvature of $M$ at the
point $x$ is strictly greater than $c.$ Using the definition of Ricci
curvatures, it follows  
\[
\theta _{k}(x)>c,\forall \ k\in \overline{2,n}. 
\]

{\bf Author's addresses} : Prof. Univ. Dr. Constantin Udri\c ste, University  Politehnica of Bucharest, Faculty of Applied Sciences,
 Splaiul Independen\c tei 313, Bucharest, 060042, Romania.

E-mail address: udriste@mathem.pub.ro

 Asistent Teodor Oprea, University of Bucharest, Faculty of Mathematics and Informatics, Str. Academiei 14,
  Bucharest, 70109, Romania.

E-mail address: teodoroprea@yahoo.com


\begin{thebibliography}{9}
\bibitem{}  R.L. Bishop, \textsl{Infinitesimal convexity implies local
convexity}, Indiana Univ. Math. J. \textbf{24} (1974), no.2, 169-172.

\bibitem{}  B.Y. Chen, \textsl{Mean curvature shape operators of isometric
immersions in real-space-forms}, Glasgow Math. J. \textbf{38 }(1996), 87-97.

\bibitem{}  B.Y. Chen, \textsl{Relations between Ricci curvature shape
operator for submanifolds with arbitrary codimensions, }Glasgow Math. J. 
\textbf{41 }(1999), 33-41.

\bibitem{}  S. Kobayashi, K. Nomizu, \textsl{Foundations of Differential
Geometry}, vol 1,2, Interscience, New York, 1963, 1969.

\bibitem{}  J.A. Thorpe, \textsl{Elementary topics in differential geometry}%
, Springer-Verlag, 1979.

\bibitem{}  C. Udri\c{s}te, \textsl{Convex hypersurfaces}, Analele \c{S}t.
Univ. Al. I. Cuza, Ia\c{s}i \textbf{32 }(1986), 85-87.

\bibitem{}  C. Udri\c{s}te, \textsl{Convex Functions and Optimization
Methods on Riemannian Manifolds}, Mathematics and Its Applications, \textbf{%
297}, Kluwer Academic Publishers Group, Dordrecht, 1994.
\end{thebibliography}
\end{document}